\input{BoxedEPS.tex}
\SetRokickiEPSFSpecial 
\HideDisplacementBoxes
\documentclass[12pt]{amsart}
\begin{document}
\title{Solvable Groups of Exponential Growth and HNN Extensions}
\author{Roger C. Alperin}
\address{\tt E-mail: alperin@mathcs.sjsu.edu,\ Department of Mathematics and Computer Science,\ 
San Jose State University,\  San Jose, CA 95192 USA
 }

\newtheorem{theorem}{Theorem}
\newtheorem{proposition}{Proposition}
\newtheorem{lemma}{Lemma}
\newtheorem{corollary}[theorem]{Corollary}
\newtheorem{conjecture}[theorem]{Conjecture}
\newtheorem{quest}[theorem]{Question}
\newtheorem{exer}[theorem]{Exercise}
\def\picture #1 by #2 (#3){\vbox to #2
      {\hrule width #1 height 0pt depth 0pt \vfill\special{picture #3}}}
\def\scaledpicture #1 by #2 (#3 scaled #4){{
  \dimen0=#1 \dimen1=#2
  \divide\dimen0 by 1000 \multiply\dimen0 by #4
  \divide\dimen1 by 1000 \multiply\dimen1 by #4
  \picture \dimen0 by \dimen1 (#3 scaled #4)}
  }

\def\cal{\mathcal}
\def\myexample#1{ \begin{exmple} {\rm #1} \end{exmple}}

\def\proof{{\bf {\noindent}Proof: }}
\newsymbol\bsq 1004
\def\endproof{\hfill${\bsq}$\bigskip}

\def\question{\par{\bf {\smallskip}{\noindent}Question: }}
\def\definition{\par{\bf {\smallskip}{\noindent}Definition: }}
\def\remark{\par{\bf {\smallskip}{\noindent}Remark: }}
\def\note{\par{\bf {\smallskip}{\noindent}Note: }}
\def\example{\par{\bf {\smallskip}{\noindent}Example: }}

\def\endex{\bigskip}
\def\endef{\bigskip}
\def\epar{\medskip}
\def\df{\em}
\def\title{\em}
\def\nindex{\index}
\def\ex{\bf}
\def\bar{\overline}

\newcommand\R{\mbox{\bf R}}
\newcommand\C{\mbox{\bf C}}
\newcommand\Z{\mbox{\bf Z}}
\newcommand\Q{\mbox{\bf Q}}

\maketitle

An extraordinary theorem of Gromov, \cite{Gv}, characterizes the finitely generated groups of polynomial growth;
a group has polynomial growth iff it is nilpotent by finite. This theorem went a long way from its
roots in the class of discrete subgroups of solvable Lie groups.  Wolf, \cite{W}, proved that a polycyclic
group of polynomial growth is nilpotent by finite. This theorem is primarily about linear groups and another
proof by Tits appears as an appendix to Gromov's paper. In fact if G is torsion free polycyclic and not nilpotent 
then Rosenblatt, \cite{Rt}, constructs a free abelian by cyclic group in G, in which the automorphism is expanding and
thereby constructs a free semigroup. The converse of this, that a finitely
generated nilpotent by finite group is of polynomial growth is relatively easy; but in fact one can also use
the nilpotent length to estimate the degree of polynomial growth as shown by Guivarc'h, \cite{Gh},  Bass, \cite{Bs}, and Wolf, \cite{W}. The
theorem of Milnor, \cite{M}, on the other hand shows that a finitely generated solvable group, not of exponential growth, is
polycyclic. Rosenblatt's version of this, \cite{Rt},  is that a finitely generated solvable group without a two
generator free subsemigroup is polycyclic.  We give another version of Milnor's theorem using the HNN construction. A
consequence is that a finitely generated solvable group $G$ has the ERF (extended residually finite)  property iff $G$ is
polycyclic.

We briefly review the HNN construction.
Generally, the HNN has a given base group $B$ and two subgroups, $H_1, H_2$ together with an (external) element $t$ of
infinite order which conjugates
$H_1$ to $H_2$, $\Gamma=<B,t\ |\ tH_1t^{-1}=H_2>.$ For solvable groups, a good example, is the
 group
$\Gamma_1=<a,t\ |\ tat^{-1}=t^2>.$ Many one relator groups have HNN decompositions; for example, consider
$\Gamma_2 = <a,t\ |\ a=[tat^{-1},t^2at^{-2}]>.$   This is, in fact, the HNN extension with base
$H=<a_0,a_1,a_2\ |\ a_0 = [a_1,a_2] >$ and free subgroups $F_1=<a_0=a,a_1=tat^{-1}>$, $F_2=<tat^{-1},
a_2=t^2at^{-2}>$ amalgamated, so that $\Gamma_2=<H ,t\ |\ tF_1t^{-1}=F_2 >$.

These HNN constructions are ascending in the sense that a conjugation of the subgroup ascends or gets strictly larger in
$\Gamma$.  We say it is ascending with base $B$ if $\Gamma$ is generated by $B$ and $t$, so that 
$$\Gamma=<B,t\ |\ tBt^{-1}\subset B>.$$
 In the first example, the
cyclic group generated by
$a$ will by repeated conjugation by
$t^{-1}$ ascend to a group isomorphic to $Z[1/2]$; the square root of $a$ exists since $(t^{-1}at)^2=a$.   In the
second example above, the subgroup
$F_1$ does not contain
$a_2$, but
$F_2$ does contain
$a_0$ so this is properly ascending with base $F_2$ and conjugation by $t^{-1}$. Since the normal subgroup generated by
$a$ is perfect, and that subgroup is locally a free group it is infinitely generated. Brown, \cite{Bn}, considers such
ascending 1-relator groups, and more general groups, in the context of actions on trees and HNN valuations.
\bigskip

It is well-known that a free product with amalgamation $A\genfrac{}{}{0pt}{}{*}{C}B$ with $C$ not of index less than or equal to 2 in each factor $A, B$,
must contain a free semigroup (and even a free group). Without loss of generality choose
$x\in A-C$, $y,z\in B-C$ distinct coset representatives, then $xy, xz$ generate a free semigroup from
the alternating word property for elements in free products with amalgamations. We do a similar
construction for ascending HNN extensions, where the base is properly included, to obtain the following.

\begin{lemma} A properly ascending HNN extension with base B,  
$$<B,t\ |\ tBt^{-1}\subset B>,$$ contains a
two generator free semigroup. 
\end{lemma}

\proof
Let $B_i=t^{i}Bt^{-i}$, $i\in Z$ and $T=\{t^n\ |\ n > 0\}$, the positive powers of $t$.  For the proof
of this lemma, we shall choose a coset representative $u\in B-B_1$, since the base is properly
ascending. Moreover, $B_1$ is normalized by $t$, but not $t^{-1}$.  Consider now the subsets
$C=B_1-\{1\}$, and $X=CT$.   Because of the  homomorphism $G\longrightarrow Z=<t>$, with kernel $K$,
the normal subgroup generated by $B$, the exponent  on $t$ is well defined when an element of $G$ is
expressed as $K<t>$. Also the subgroup $K$ is a properly ascending union of subgroups, $B_i\subset
B_{i-1}$, for $i\in Z$. 

The translate $tX$  is a subset of $X$, since for $c\in C$, then $tct^{-1}\in C$ so that
$t(ct^i)=(tct^{-1})t^{i+1}$. Consider also $tuX$; the element
$tuct^i$=$(tut^{-1})(tct^{-1})t^{i+1}$ is contained in $X$; moreover since the exponent on $t$ is
positive and  $G$ is properly ascending, for $u\in
B-B_1$, the element $(tut^{-1})(tct^{-1})$ is not 1.  Finally,  $tX$ and $tuX$   are
disjoint subsets since if $tu(ct^i)=t(dt^j)$ $c,d\in C$ then $uc=d$, which is impossible  by the
choice of $u$.

Now from the disjointness of $X_1=tX$ and $X_2=tuX$, and also $X_1\cup X_2\subset X$, it now easily
follows that $t, tu$ generate a free semigroup. Any distinct words $w_1, w_2$ in  $t, tu$ without
loss of generality begin on the left with $t, tu$ respectively, so $w_1X\subset X_1, w_2X\subset X_2$
are different.\endproof

\medskip
Since we can add free semigroup generators to a generating set the following is immediate.

\begin{proposition} A group which contains a free semigroup on two generators has exponential
growth.
\end{proposition}
\bigskip
Our version of Milnor's Theorem is the following.

\begin{theorem}
\label{Theorem A}  A finitely generated solvable group which is not polycyclic contains a subgroup
which is a properly ascending HNN extension.
\end{theorem}
\proof Suppose that $\Gamma$ is a finitely generated solvable group. Descending down the
solvable series we obtain finitely generated layers up to level $n$, say, and infinitely generated mod
$\Gamma^{n+2}$ at $n+1$ or else the group is polycyclic. Consider the finitely generated solvable
group $G=\Gamma^n/\Gamma^{n+2}$;  $P=\Gamma^n/\Gamma^{n+1}$ is finitely generated abelian, and 
$A=\Gamma^{n+1}/\Gamma^{n+2}$ is infinitely generated abelian. Also $P$
is not finite, since otherwise $A$ is finitely generated. We may assume that $P$ is torsion free
by passing to a subgroup of finite index in $G$. It suffices now to prove that
$G$ contains a properly ascending HNN extension. 

We claim that there is a $t\in G$ of infinite order, and subgroup 
$B\subset A$,  so that $tBt^{-1}\subset B$ is proper. Let $t$ be an arbitrary element of $G$ of
infinite order which maps non-trivially to $P$. Let $R=Z[t]$ be the group ring of the monoid
generated by $t$. Let $a$ be an arbitrary element of $A$; and let $M=Ra\subset A$ be the cyclic module
generated by $a$. Certainly using `module' action, $tM\subset M$; if this inclusion is proper we are
done, $B=M$. So we may assume now that for every element $t$ of $P$, $tM=M$; but then also $t^{-1}M=M$. Thus, there 
is some polynomial in $t^{-1}$ with constant term equal to one which annihilates $a$. Hence a monic polynomial in $t$
annihilates $a$. It follows that $M$ is finitely generated as an abelian group. Thus
using $t$ and $t^{-1}$ as above we have that for any element
$a$, the set of all conjugates $t^iat^{-i}, i\in Z$ is a finitely generated abelian subgroup of $A$.

Notice that since $P$ is finitely presented, and $G$ is finitely generated, then $A$
is finitely generated as a normal subgroup of $G$. Let $\{a_1, a_2, ... , a_m\}$ be the finite set of
normal generators of $A$; now using the finite set of generators $\{t_1,t_2,...,t_n\}$ of $G$ mod $A$,
we obtain a finite set of generators for $A$. Consider $\{t_1^ka_it_1^{-k}:i=1,...m,  k\in
Z\}$, since we have already shown that for any $a \in A$, $t^iat^{-i}, i\in Z$ is finitely
generated, this set is generated by a finite set of elements. Now continue with $t_2$ up to 
$t_n$. This gives a finite set of generators for $A$ which contradicts our assumption on $A$; thus
there must be a properly ascending HNN extension as a subgroup of $G$. \endproof

\medskip
We extract the following Burnsidesque local-global property  of a  group $\Gamma$: 
\par\noindent {\bf Property
$\pi$}: Given any finitely generated module $A$ for $\Gamma$; if for every $t\in\Gamma$, $a\in A$, the
$Z[t]$ module generated by $a$ is finitely generated then $A$ is a finitely generated abelian group. 
\medskip
Polycyclic groups
have property $\pi$. If $\Gamma$ is polycyclic we can use the normal form of elements, given from 
the polycyclic decomposition, to build iteratively as above, a finite set of generators
starting from the finite set of module generators. In a similar way, we see that any group which has the following {\it
bounded generation} property, has property $\pi$: there are a finite set of elements $F$, and a fixed integer $N$, so that
every element of $\Gamma$ has an expression as
$g_1^{m_1}g_2^{m_2}...g_N^{m_N}$, $g_i\in F$, $m_i\ge 1$. Many lattices in semisimple Lie groups have this bounded
generation property; it is related to a positive solution to the congruence subgroup problem, \cite{Rk}. Also, it is immediate,
that a finitely presented torsion group with property $\pi$ is finite. At the conference, J. S. Wilson kindly pointed out to
me the following references. Kropholler \cite{K} has shown that finitely generated minimax solvable groups have bounded generation.
This includes the Baumslag-Solitar group
$\Gamma_1$ discussed earlier. Brookes \cite{Br}, has investigated Engel elements of solvable groups and its relation to property
$\pi$.
\bigskip
The following is immediate from the proof and discussion given above.
\begin{theorem}
\label{Theorem B}   A finitely generated group $G$ containing a normal subgroup $N$,
with infinitely generated abelianization $N/N^{'}$, contains a properly ascending HNN extension if
 $G/N$ is a  finitely presented group having property $\pi$. 
\end{theorem}

In particular this applies to $G=F/R^{'}$ for an infinite group with finite (free) presentation $F/R$ having
property $\pi$. One can also replace $R^{'}$ by any other term $R^{(n)}, n\ge 2$ in the derived series of $R$ to obtain similar
results.
\medskip

 Recall that the group $G$ has the property ERF
 if every subgroup is closed in the profinite topology, or equivalently, given any subgroup $S$ and any
element $x\in G-S$, there is a finite index subgroup (equivalently normal subgroup) containing $S$ and not $x$; this is also
the same as the existence of a finite quotient of $G$ so that $S$ is represented trivially and  $x$ non-trivially. 
The property LERF requires that only finitely generated subgroups are closed in the profinite topology. 
It has been shown by Malcev, \cite{Mv}, that a group which is a split extension,  with normal subgroup having the ERF property, and
quotient which has LERF, itself is ERF. This shows a large class of solvable groups have the property ERF. Since polycyclic
groups satisfy the maximal condition on subgroups, all subgroups are finitely generated, we see immediately that a polycyclic
group has the property ERF.

One can also apply our main Theorem \ref{Theorem A} to obtain the following surprising result.

\begin{theorem}
\label{Theorem C}   A finitely generated solvable group $G$ has the property ERF iff $G$ is polycyclic.
\end{theorem}

\proof  As remarked above, polycyclic groups have the property ERF. On the other hand, if
the group is not polycyclic, we can find a subgroup
$B$ and element $t$ so that $tBt^{-1}\subset B$, and $B_1=tBt^{-1}\neq B$.  Therefore in a finite image an element
$a\in B-B_1$ can not be separated from $B_1$.
\endproof

In particular, it follows from this (or the maximal condition) that a polycyclic group is never a properly ascending HNN
extension.

\end{document}